\documentclass[12pt]{article}
\usepackage[colorlinks=true]{hyperref} 
\usepackage{url}      
\usepackage{doi}      
\usepackage{amsmath, amsthm} 
\usepackage{amsfonts} 
\usepackage{amssymb} 
\usepackage{dirtytalk} 
\usepackage{esint} 
\usepackage{subcaption} 
\usepackage[colorlinks=true]{hyperref} 
\usepackage{tikz} 
\usepackage{mathtools} 
\usepackage{enumerate} 
\usepackage{bm} 
\usepackage{authblk} 
\usepackage[nameinlink,noabbrev,capitalize]{cleveref} 
\crefname{ass}{Assumption}{Assumptions}
\usepackage[left=2.54cm,right=2.54cm,top=2.54cm,bottom=2.54cm]{geometry} 
\usepackage{comment} 

\newcommand{\C}[1]{\mathcal{#1}}

\newcommand{\B}[1]{{\bf #1}}

\newcommand{\ad}{{\rm{ad}}}
\newcommand{\mt}{{-i}}

\newcommand{\Fx}{{\rm Fix}(\B{X})}
\newcommand{\expl}[1]{{ \color{blue} }}

\usepackage{pgfplots} 
    \pgfplotsset{compat=newest}
    \usetikzlibrary{plotmarks} 
    \usetikzlibrary{calc} 
    \usetikzlibrary{positioning} 
    \usetikzlibrary{arrows.meta} 
    \usetikzlibrary{matrix} 
    \usepgfplotslibrary{patchplots} 
    \usetikzlibrary{shapes.geometric} 
\numberwithin{equation}{section}

\newcommand*\EX{\hspace*{0em plus 1fill}\makebox{\qedsymbol}}

\newtheorem{theorem}{Theorem}[section]
\newtheorem{prop}[theorem]{Proposition}
\newtheorem{ass}[theorem]{Assumption}
\newtheorem{definition}[theorem]{Definition}

\newtheorem{remark}[theorem]{Remark}

\theoremstyle{definition}
\newtheorem{example}[theorem]{Example}
\title{Structure versus regularity of set-valued maps in convex generalized Nash equilibrium problems in Banach spaces\thanks{Both authors acknowledge support by the Deutsche Forschungsgemeinschaft (DFG), SFB-TRR 154, project ID 239904186, subproject C08. MH also acknowledges support by MATH+: The Berlin Mathematics Research Center, Excellence Cluster EXC-2046/1, project ID 390685689, subproject AA4-13.}}
\author[1]{Marcelo Bongarti\footnote{Corresponding author: Marcelo Bongarti, bongarti@wias-berlin.de}}
\author[1,2]{Michael Hintermüller}
\affil[1]{Weierstraß Institute for Applied Analysis and Stochastics, Berlin, Germany}
\affil[2]{Humboldt-Universität zu Berlin, Berlin, Germany}
\date{
}

\begin{document}
\maketitle

\begin{abstract}
A generalized Nash equilibrium problem (GNEP) in Banach space consists of $N > 1$ optimal control problems with couplings in both the objective functions and, most importantly, in the feasible sets. We address the existence of equilibria for convex GNEPs in Banach space. We show that the standard assumption of lower semicontinuity of the set-valued constraint maps – foundational in the current literature on GNEPs – can be replaced by graph convexity or the so-called Knaster–Kuratowski–Mazurkiewicz (KKM) property. Lower semicontinuity is often essential for obtaining upper semicontinuity of best response maps, crucial for the existence theory based on Kakutani–Fan fixed-point arguments. However, in function spaces or in settings with partial differential equation (PDE) constraints, verifying lower semicontinuity becomes much more challenging (even in convex cases), whereas graph convexity, for example, is often straightforward to check. Our results unify several existence theorems in the literature and clarify the structural role of constraint maps. We also extend Rosen’s uniqueness condition to Banach spaces using a multiplier bias framework.
\end{abstract}

\section{Introduction}

Game theory provides a mathematical framework for studying interactions among decision makers. Games can be cooperative, if enforceable agreements are allowed, or noncooperative \cite{harsanyi_general_1988}, and since the seminal work of von Neumann and Morgenstern \cite{von_neumann_theory_2007} and Nash's foundational equilibrium concepts \cite{nash_equilibrium_1950,nash_non-cooperative_1951}, equilibrium analysis has become a central paradigm in economics, engineering, and the social sciences. A \emph{generalized Nash equilibrium problem} (GNEP) extends the classical setting by allowing each player's feasible set to depend on the strategies of others. This formulation captures shared or coupled constraints, which arise naturally in models such as electricity markets \cite{hobbs_nash-cournot_2007,jing-yuan_spatial_1999}, 
environmental management \cite{breton_game-theoretic_2006} and resource allocation \cite{roughgarden_algorithmic_2010}. 

Most classical results on GNEPs are finite-dimensional. However, important applications such as the coupling of energy markets to physical transport of energy carriers (see, e.g., \cite{bongarti_optimal_2024} and the references therein) lead to GNEPs with partial differential equations (PDEs) modeling the transport. Such infinite-dimensional formulations, particularly those constrained by PDEs, remain less developed due to compactness and regularity issues. While some works analyze PDE-constrained Nash equilibrium problems via optimality systems \cite{roubicek_nash_2007,roubicek_noncooperative_1999,ramos_nash_2002,ramos_pointwise_2002,MR3032892}, genuine generalized frameworks for multiple players have been addressed mainly under jointly convex assumptions \cite{hintermuller_generalized_2015,hintermuller_generalized_2024,gahururu_generalized_2022,gahururu_risk-neutral_2023,hintermueller_pde-constrained_2013,kanzow_multiplier-penalty_2019}. Beyond standard compactness of the feasible strategy sets and closedness of constraint maps, these studies typically require lower semicontinuity of the constraint maps to invoke Kakutani–Fan-type fixed-point theorems. While this is a standard assumption in abstract analysis, it is often not verifiable in function space settings.

This paper contributes to the existence theory for convex GNEPs in Banach spaces by clarifying the trade-off between structural assumptions, such as convexity and regularity of set-valued constraint maps. We show that equilibrium existence can be ensured under \emph{graph-convexity} or the \emph{Knaster}–\emph{Kuratowski}-\emph{Mazurkiewicz (KKM)} property without assuming lower semicontinuity. These geometric conditions are often much easier to verify in infinite-dimensional settings, in particular in PDE-constrained contexts and thus broaden the class of problems where equilibrium existence can be rigorously established. This relaxes one of the most restrictive assumptions in the standard theory while retaining the classical framework as a special case.

We also show that Rosen’s uniqueness condition for variational equilibria \cite{rosen_existence_1965} extends naturally to Banach spaces, and introduce the notion of \emph{multiplier bias} as a mechanism for equilibrium selection. 

Compared with the existing literature on convex GNEPs and quasi-variational inequalities (QVIs) \cite{facchinei_generalized_2007,harker_generalized_1991,pang_quasi-variational_2005,chan_generalized_1982}, the novelty is not in a single new fixed-point theorem, but in the way how structural assumptions on the constraint maps and the objective functionals can be traded against regularity assumptions, with an eye on PDE-constrained applications. The paper thus complements, rather than overlaps, with the jointly convex theory in \cite{hintermuller_generalized_2015,kanzow_multiplier-penalty_2019} and with equilibrium results in function spaces \cite{mas-colell_equilibrium_1991}.

The rest of the paper is organized as follows:
\Cref{gen_gam} introduces the mathematical structure of generalized games in Banach spaces. 
\Cref{convex_reg} establishes equilibrium existence under graph-convexity and KKM assumptions. \Cref{bias} extends Rosen’s uniqueness result and interprets it as a selection mechanism through multiplier bias. 

\section{Generalized finite games in Banach spaces} \label{gen_gam}

Let $N>1$ be a fixed natural number and $I = \{1, \cdots, N\}$ be a set of indices representing the players or agents participating in the game. 

The strategy set of each player is a fixed Banach space $X_i$, likely of infinite dimension, with its topological dual denoted by $X_i^*$. Representing their private (or individual) constraints, each player observes a nonempty, closed and convex subset $X_i^\ad \subset X_i$. We write
$X \coloneqq X_1 \times \cdots \times X_N$, analogously for $X^*$,
and $X^\ad \coloneqq X_1^\ad \times \cdots \times X_N^\ad$
to denote the full strategy space, its topological dual, and the admissible strategy set, respectively. 

In game theory, it is common to use the notation $x_\mt$ to denote the bundle of $N-1$ strategies containing all strategies of a given bundle $x \in X$ \emph{except} that of player $i.$ Moreover, we write $x = (x_i, x_\mt)$ when emphasis on player $i$'s strategy is needed, but without changing the original order of the components of $x$. Consequently, we let $X_\mt \coloneqq X_1 \times \cdots \times X_{i-1} \times X_{i+1} \times \cdots \times X_N$ (and the corresponding $X_\mt^\ad$) and say that $x_\mt \in X_\mt.$ In earlier literature on Nash games the bundle $x_\mt$ is sometimes referred to as an \emph{$i$-incomplete combination} \cite{harsanyi_general_1988}.

To complete the description of a generalized game, two ingredients are still needed. First, a family of objectives $\C{J} = \{\C{J}_i\}_{i = 1}^N$ with each function $\C{J}_i$ being a real-valued function defined on the open set $O \coloneqq O_1 \times \cdots \times O_N$ where each $O_i$ is an open set of $X_i$ containing $X_i^\ad.$ And second, a family of set-valued constraint maps $\{\B{X}_i\}_{i=1}^N$ with each $\B{X}_i$ defined from $O_\mt$ to $O_i$. Then $\B{X}_i(x_\mt)$ is the set of \emph{admissible} strategies which player $i$ can take given the other players decisions $x_\mt$. Finally, we denote by $\B{X}$ the set-valued map $\B{X}_1 \times \cdots \times \B{X}_N.$

\begin{definition}[\bf GNEP]
	A generalized Nash equilibrium problem (GNEP), denoted by $G = (\B{X}, \C{J})$ is a family of $N$ coupled constrained optimization problems of the form
	\begin{equation}\label{gamei}\tag{$P_i$}
		\begin{cases}
			\text{Given } x_\mt \in X_\mt^\ad; \\
			\text{Minimize } \C{J}_i(x_i, x_\mt) \ \text{over } x_i \in X_i^\ad \cap \B{X}_i(x_\mt). 
		\end{cases}
	\end{equation}
\end{definition}

If $\B{X}_i(x_\mt) = O_i$ for all $x_\mt \in X_\mt$ and all $i \in I$, then a GNEP reduces to a so called Nash equilibrium problem (NEP) as dependence on $x_{-i}$ only occurs in the objectives, but not the constraints.

For future reference, we denote by $\C{D}_i$ the \emph{domain} of the map $\B{X}_i$, i.e., 
\begin{equation}\label{domain} \C{D}_i \coloneqq \{x_\mt \in X_\mt^\ad; \ \B{X}_i(x_\mt) \neq \emptyset\}.
\end{equation}

As a first example, the simplest situation in which generalized games appear is when a constraint is \emph{shared} among players. This is the case when, for example, a \textit{regulatory} individual or agency determines that the strategy bundle $x$ must belong to a certain \emph{closed} set $\C{C} \subset X.$ This leads to non-constant constraint maps 
\begin{equation}
	\label{constraint_map}
	\B{X}_i(x_\mt) = \{x_i \in X_i; x = (x_i, x_\mt) \in \C{C}\} = \pi_i(\C{C}),
\end{equation}
where $\pi_i$ denotes the projection onto $X_i$. 

We are interested in the concept of \textit{Nash equilibrium}  \cite{nash_equilibrium_1950,nash_non-cooperative_1951} as the notion of \emph{solution} for a generalized game.

\begin{definition}[\bf GNE] 
	We say that $\overline{x} \in \B{X}(\overline{x})$ is a generalized Nash equilibrium (GNE) of a game $G = (\B{X}, \C{J})$ provided
	\begin{equation}\label{analytic_min_pn}
		\C{J}_i(\overline{x}) = \min\limits_{x_i \in X_i^\ad \cap \B{X}_i(\overline{x}_\mt)} \C{J}_i(x_i, \overline x_i)
	\end{equation}
	for all $i \in I.$ We denote the set of generalized Nash equilibria by $\B{E}(G).$
\end{definition}

Very simple finite dimensional examples can illustrate the fact that constant versus non-constant constraint maps may drastically change the set of equilibria, see e.g. \cite{facchinei_generalized_2007}.
Also, the issue of compacteness (of the feasible strategy sets) becomes essential in games with nonsmooth objectives \cite{mas-colell_equilibrium_1991}.
Furthermore, in infinite dimensions, even compactness cannot ensure that an equilibrium exists, see \cite{mas-colell_equilibrium_1991} for examples. 

In particular for PDE-constrained games, the natural topology in which state equations are well posed, closed and bounded sets are usually only weak or weak-star compact, hence lower semicontinuity of the induced constraint maps with respect to that topology is often impractical. This motivates the move to structural conditions on the graph, as we propose in this paper.

\begin{definition}[\bf Convex GNEP] We say that the generalized game $G = (\B{X}, \C{J})$ is convex if, for all $i \in I$, the functions $\B{X}_i$ are convex-valued and the functions $\C{J}_i(\cdot, x_\mt)$ are convex.
\end{definition}

Let $X,Y$ be Banach spaces. The graph of a set-valued map $F: X \to 2^Y$ is defined as $\B{Gr}(F) = \{(x,y) \in X \times Y: y \in F(x)\}$. We say that $F$ is lower semicontinuous (l.s.c) at $x_0$ if for any open subset $V$ of $Y$ such that $F(x_0) \cap V \neq \emptyset$ there exists an open set $U$ of $X$ containing $x_0$ and such that $F(x) \cap V \neq \emptyset$ for all $x \in U.$ Equivalently, $F$ is l.s.c. at $x_0$ if for any $y_0 \in F(x_0)$ and any sequence $x_n$ in $X$ converging to $x_0$ there exists a sequence $y_n$ in $Y$ converging to $y_0$ and such that $y_n \in F(x_n).$ $F$ is said to be l.s.c on a subset $\C{X}$ of $X$ if it is l.s.c at all $x_0 \in \C{X}.$ For more details, see \cite{beer_topologies_1993}.

Along with lower semicontinuity of the constraint maps $\B{X}_i$ (for each $i$), the following two assumptions are standard in the study of GNEPs.
\begin{ass}\label{ass_compactness_regularity} 
	For each $i \in I$, 
	\begin{itemize} 
		\item[\bf (i)] $X_i^\ad$ is nonempty, compact and convex;
		\item[\bf (ii)] the map $\B{X}_i: X_\mt \to 2^{X_i}$ is convex-valued and has closed graph;
		\item[\bf (ii')] $\B{X}_i(x_\mt) \subset X_i^\ad$ for all $x_\mt \in X_\mt^\ad.$
	\end{itemize}
\end{ass}
\begin{ass}\label{ass_objectives} 
	For each $i \in I$ and $x_\mt \in X_\mt$
	\begin{itemize}
		\item[\bf (i)] the function $\C{J}_i: X \to \mathbb{R}$ is continuous; 
		\item[\bf (ii)] the function $\C{J}_i(\cdot, x_\mt): X_i \to \mathbb{R}$ is convex.
	\end{itemize}
\end{ass}
\Cref{ass_compactness_regularity} above plays a pivotal role in the study of convex GNEPs. Established methodologies, such as the K. Fan inequalities \cite{fan_minimax_1953} and various formulations of Kakutani's theorem \cite{beer_topologies_1993}, rely heavily on this assumption. It is worth mentioning that (ii') is not actually necessary for any of the results discussed here, but we make it just so there is no need to write $X_i^\ad \cap \B{X}_i(x_\mt)$ everywhere. 

On the other hand, \Cref{ass_objectives} fixes the structure of the objective functions considered in this paper. See, for example, \cite{arrow_existence_1954,mas-colell_equilibrium_1991} for a justification of these conditions in the context of economics.

Assuming lower semicontinuity of $\B{X}_i$ is crucial for the closedness of the so called graph of the best-response map. Fixed points of such maps characterize Nash equilibria and this is why Kakutani's theorem is central in the current existence theory for GNEPs.
This assumption is, however, rarely easy to verify in practical settings such as PDE-constrained or function space games, see \cite{bensoussan_impulse_1987,hintermuller_recent_2019,facchinei_generalized_2007,facchinei_decomposition_2011}.  
 This limitation partly explains the popularity of models where the coupling constraints have the structured form \eqref{constraint_map}, for which the lower semicontinuity issue can be bypassed by exploiting properties of the Nikaido-Isoda function. This shows that to a certain extent, regularity (lower semicontinuity) can be traded by structure (jointly convex games). 

Our main result in this paper extends this idea to a general class of convex GNEPs (not only jointly convex). It says that if, in addition to convexity of GNEP, the graph of the bundle constraint map $\B{X}$, i.e. the set $\{(x,y) \in X \times X: y \in \B{X}(x)\}$ is convex, then {\it lower semicontinuity of the maps $\B{X}_i$ is no longer necessary for existence of an equilibrium.} Details on this are discussed in the next section.

\section{Potential lack of lower semicontinuity of constraint maps}\label{convex_reg}

For convex games we show next that lower semicontinuity of the constraint maps is not needed for equilibrium existence once additional geometric structure is available. This allows us to cover games with significantly more complex and heterogeneous constraints, such as generalized PDE-constrained games where each player’s state equation depends on forecasts of incomplete or partially available information.

\begin{definition}
	Let $X$ be a nonempty subset of a topological vector space $L$. We say that a set valued map $F: X \to 2^L$ is graph-convex if its graph is convex.
\end{definition}

\begin{theorem}\label{no_lowersemicontinuity_theo_0}
	Assume that $\B{X}$ is graph-convex and that $\Fx$, the set of feasible fixed points of $\B{X}$, i.e., $$\Fx \coloneqq \{x \in X^\ad; x \in \B{X}(x)\}$$ is nonempty. Then the game $G = (\B{X},\C{J})$ has an equilibrium.
\end{theorem}

It is important to emphasize here that graph-convexity does not \emph{replace} lower semicontinuity as a weaker condition in general. In fact, it is well known that graph-convexity implies lower semicontinuity in the interior of the domain. A proof of this fact can be found in \cite[Theorem 5.9(b)]{rockafellar_variational_1998} in finite dimensions, and the same proof can be adapted to the Banach space case without much work. However, in infinite dimensions interiority is delicate, as we illustrate in the example below.

\begin{example}
	Let $I=\{1,2\}$ and $X_1 = X_2 = L^2(0,1)$. Set $X := X_1\times X_2$ and define the closed subspace
	$$
	C := \left\{\,u\in L^2(0,1) : \int_0^1 u(x)\,dx = 0\,\right\}.
	$$
	Since $C = \ker \ell$ for the continuous, surjective linear functional
	$$
	\ell(u) := \int_0^1 u(x)\,dx,
	$$
	the space $C$ has codimension~$1$ in $L^2(0,1)$ and hence empty interior in the strong $L^2$ topology.
	
	Define the closed convex set
	$$
	K := \left\{v\in L^2(0,1): v \geqslant 0 \ \mbox{a.e.}, \ \|v\|_{L^2}\leqslant 1\right\},
	$$
	where \say{a.e.} stands for \say{almost everywhere} (in the sense of Lebesgue), and fixed feasible sets
	$$
	X_i^\ad = X_i = L^2(0,1),\qquad i=1,2.
	$$
	Finally, the player-specific constraint maps $\B{X}_i : X_\mt \to 2^{X_i}$ are given by
	$$
	\B{X}_i(x_\mt) :=
	\begin{cases}
		K, & \text{if }\displaystyle \int_0^1 x_\mt(x)\,dx = 0,\\[3pt]
		\emptyset, & \text{if }\displaystyle \int_0^1 x_\mt(x)\,dx \neq 0,
	\end{cases}
	$$
	
	Thus each player's feasible set depends on the other player's strategy. The joint constraint map
	$
	\B{X}(x_1,x_2) := X_1(x_1,x_2)\times X_2(x_1,x_2)
	$
	is then given by
	$$
	\B{X}(x_1,x_2) =
	\begin{cases}
		K_1\times K_2, &\text{if } x_1\in C,\ x_2\in C,\\[3pt]
		\emptyset, & \text{otherwise}.
	\end{cases}
	$$
	(One may take, for instance, $\C{J}_i \equiv 0$ such that the objectives play no role in this example.)
	
	The graph of the joint constraint map is
	$$
	\B{Gr}(\B{X})
	= (C\times C)\times (K_1\times K_2),
	$$
	which is convex because $C$, $K_1$, and $K_2$ are convex, respectively. The domain of $\B{X}$ is
	$$
	{\rm dom} \ \B{X} = \left\{\,x\in X : \B{X}(x)\neq\emptyset\,\right\}
	= C\times C,
	$$
	which has empty interior in $X$.
	
	Further, lower semicontinuity of $\B{X}$ fails at every $\overline{x}
    \in{\rm dom } \ \B{X}$. In fact, fix any $\overline{y}\in \B{X}(\overline{x}) = K_1\times K_2$. Since $C$ is a proper closed subspace of $L^2(0,1)$, we can choose sequences
	$$
	x_1^n \to \overline{x}_1 \quad\text{with } x_1^n\notin C,\qquad
	x_2^n \to \overline{x}_2 \quad\text{with } x_2^n\notin C.
	$$
	Set $x^n:=(x_1^n,x_2^n)$. Then $x^n\to\overline{x}$ in $X$, but by construction
	$
	\B{X}(x^n) = \emptyset$ for all $n$.
	Hence there is no sequence $y^n\in\B{X}(x^n)$ converging to $\overline{y}$, i.e., lower semicontinuity of $\B{X}$ fails at $\overline{x}$.
	\EX
\end{example}
	The example above is a genuine two-player generalized Nash game whose joint constraint map has a convex graph, but a domain with empty interior in the strong topology, and it illustrates how lower semicontinuity can fail in such \say{thin} situations. 

    \subsection{The Nikaido-Isoda Function}

The Nikaido-Isoda function $\Psi: X \times X \to \mathbb{R}$ reads 
\begin{equation}\label{NI_Function} 
	\Psi(x,y) = \sum\limits_{i = 1}^N \C{J}_i(y_i, x_\mt).
\end{equation} 
Under \Cref{ass_objectives}, $\Psi$ is continuous and, for each $x \in X$, the map $\Psi(x, \cdot)$ is convex. A strategy bundle $\overline{x} \in X^\ad$ is a Nash equilibrium of the game $G = (\B{X}, \C{J})$ if and only if $\overline{x} \in \B{X}(\overline{x})$ and $\Psi(\overline x, \overline x) \leqslant \Psi(\overline x, y)$, for all  $y \in \B{X}(\overline x)$. The forward implication of this equivalence is trivial. For the backwards one, assume that $\overline{x}$ is a fixed point of $\B{X}$ and that $\Psi(\overline x, \cdot)$ reaches its global minimum at $\overline x.$ If $\overline{x}$ is not a Nash equilibrium, then for some $i \in I$, there exists $y_i \in \B{X}_i(\overline{x}_\mt)$ such that 
	$   \C{J}_i(\overline{x}) > \C{J}_i(y_i, \overline{x}_\mt)$.
	Noticing that $(y_i, \overline{x}_\mt) \in \B{X}(\overline x)$, the result follows because 
	\begin{equation*}
		\Psi(\overline x, \overline x) - \Psi(\overline x, (y_i, \overline{x}_\mt)) = \C{J}_i(\overline{x}) - \C{J}_i(y_i, \overline{x}_\mt) > 0,
	\end{equation*}
	which contradicts minimality of $\overline{x}.$

In other words, $\overline{x} \in \B{X}(\overline x)$ is a GNE if and only if $\Psi(\overline x, \overline x) = \Phi(\overline x)$ where $\Phi$ is defined as
\begin{equation}
	\label{min_prob_eq} \Phi(x) =  \min \left\{\Psi(x, y) \ \ {\rm s.t. } \ \ y \in \B{X}(x)\right\}.
\end{equation} 

Among other consequences, this equivalence allows one to characterize Nash equilibria -- under extra regularity of the objective functions -- as the solution of a quasi-variational inequality (QVI). 
It also implies that a bundle $\overline{x}$ is a GNE for a game if and only if 
\begin{equation}
	\label{fix_def_nash} 
	\overline{x} \in \left\{\hat y \in \B{X}( \overline{x}) \cap \Fx: \Psi( \overline{x},\hat y) = \Phi( \overline{x})\right\}.
\end{equation} 
We can then prove \Cref{no_lowersemicontinuity_theo_0} in a rather elementary way; it extends the classical argument of Nikaido and Isoda \cite{nikaido_note_1955} to the present setting.

\subsection{Proof of \Cref{no_lowersemicontinuity_theo_0}}
	Suppose that a GNE does not exist. 
    This means, in particular, that for each $x \in \Fx$ there exists $y_x \in \B{X}(x) \cap \Fx$ such that $\Psi(x, x) > \Psi(x,y_x).$ Fixing $y_x$, it follows from continuity of $\Psi$ that the set 
	$
	\Psi_{x} \coloneqq \{z \in X:\,\Psi(z, z) > \Psi(z,y_x)\}
	$
	is nonempty (in fact it contains $x$) and open in $X$. 
	
	Now, notice that the closedness of the graph of each map $\B{X}_i$ (see \Cref{ass_compactness_regularity}(b)) implies that the set $\Fx$ is closed. In fact, let $(x^n)$ be a sequence in $\Fx$ such that $x^n \to x \in X.$ We have $x_i^n \in \B{X}_i(x_\mt^n)$, i.e., $x^n \in {\rm Gr}(\B{X}_i)$ for each $n$ and each $i$. By closedness of the graph, it follows that $x \in {\rm Gr}(\B{X}_i)$ for each $i$, hence $x \in \Fx.$ Therefore, as a closed subset of $X^\ad$, which is compact, $\Fx$ is also compact. Hence, the trivial inclusion 
	$$
	\Fx \subset \bigcup_{x \in \Fx} \Psi_x
	$$
	implies that there exist $x_1, \cdots, x_r \in \Fx$ (with $r \in \mathbb{N}$) such that $\Fx \subset \bigcup_{i=1}^r \Psi_{x_i}$. This implies, in particular, that 
	$$
	\Psi(x,x) > \min_{1 \leqslant  i \leqslant  r} \Psi(x, y_{x_i}) \qquad \text{for all } x \in \Fx.
	$$
	
	Define, for each $1 \leqslant i \leqslant r$, the (continuous) function $g_i: X \to \mathbb{R}$ as 
	$$
	g_i(x) = -\min\{\Psi(x,y_{x_i})-\Psi(x,x), 0\}
	$$
	and notice that $g_i(x) \geqslant 0$ with $\sum_i g_i(x) > 0$ for all $x \in \Fx$. Hence the function $G$ defined by 
	$$
	G(x) = \dfrac{\sum_{i=1}^r g_i(x)y_{x_i}}{\sum_{i=1}^r g_i(x)}
	$$
	maps $\Fx$ into $A = {\rm co}(\{y_{x_1}, \cdots, y_{x_r}\})$, the convex hull of the set $\{y_{x_1}, \cdots, y_{x_r}\}$, and, in particular, $A$ to itself.  The set $\Fx$ is convex and contains each $y_{x_i}$, so the convex hull $A$ is a nonempty convex subset of $\Fx$. Since the $y_{x_i}$ span a finite-dimensional subspace of $X$, the set $A$ is compact in the induced topology. The map $G$ is continuous on $\Fx$ (and hence on $A$) because it is built from finitely many continuous functions $g_i$ and finitely many fixed points $y_{x_i}$. Therefore $G:A\to A$ is a continuous self-map of a nonempty compact convex subset of a finite-dimensional space. Thus, Brouwer's fixed point theorem applies. It follows that $G$ has a fixed point $\hat x \in A$. For indices $i$ such that $g_i(\hat x) > 0$ we have $\Psi(\hat x, y_{x_i}) < \Psi(\hat x, \hat x)$ and this, combined with convexity of $\Psi(\hat x,\cdot)$, yields 
	\begin{align*}
	\Psi(\hat x, \hat x) 
	&= \Psi\Big(\hat x, G(\hat x)\Big)
	\\ &\leqslant  \sum_{i=1}^r \frac{g_i(\hat x)}{\sum_{j=1}^r g_j(\hat x)} \,\Psi(\hat x,y_{x_i})
	< \Psi(\hat x,\hat x),
	\end{align*}
	which is a contradiction. Hence a Nash equilibrium must exist. \EX

The example below, inspired by the viscosity-regularized spot market system of \cite[Section 5.2]{hintermuller_generalized_2015}, illustrates how equilibrium existence based on graph-convexity
applies even when the feasible sets are not jointly convex.

\begin{example}[\bf PDE-constrained spot market]
	Let $\Omega = (0,1)$, $T>0$, and 
	$Q:=\Omega\times(0,T)$.
	For each player $i\in I=\{1,\dots,N\}$, the decision variable is
	$u_i\in X_i:=L^2(0,T;L^2(\Omega))$ with simple box constraints
	$0\leqslant  u_i(x,t)\leqslant \bar u_i$ a.e.\ in $Q$ for given $\bar{u}_i\in\mathbb{R}_+$.
	Each strategy bundle $u = (u_1, \cdots, u_N)$ produces a state $y\in Y:=L^2(0,T;H_0^1(\Omega))$ that satisfies the linear PDE
	\begin{equation}\label{mh:heat}
		\begin{cases}y_t - \Delta y = \sum_{i=1}^N u_i \quad\text{in } Q, \\ 
		y=0 \text{ on }\{0,1\}\times(0,T),
		y(\cdot,0)=0.
        \end{cases}
	\end{equation}
	The solution operator $S:X\to Y$ associated with $\eqref{mh:heat}$ is linear and continuous. Note that above $L^2(\Omega)$, $H_0^1(\Omega)$ denote the usual Lebesgue and Sobolev spaces \cite{MR2424078} and $L^2(0,T; L^2(\Omega))$ as well as $L^2(0,T;H_0^1(\Omega))$ are the associated Bochner spaces \cite{MR895589}. 
	
	Each player $i$ enforces the shared constraint set by an individual condition that can be interpreted as a personal buffer or uncertainty allowance.
	For example, assuming player $i$'s state requirement set (given a competition bundle $u_\mt$) is the closed convex set $\mathcal{K}_i(u_\mt)$,
	we define his or her feasible set as the convex set
	$$ \B{X}_i(u_{-i}) :=\left\{u_i\in X_i^\ad : y = S(u) \in \C{K}_i(u_\mt)\right\}. $$
	
	Hence the product mapping $\B{X}:=\B{X}_1\times\cdots\times\B{X}_N$
	is graph-convex on $X^\ad$. Hence, as long as $\Fx \neq \emptyset$ (which is a standard assumption to make) one can easily prove the existence of an equilibrium as long as the objective functions satisfy \Cref{ass_objectives}. \EX
\end{example}

If $\B{X}$ is not graph-convex, but is a KKM map, then Nash equilibria also exist.
\begin{definition}
	Let $X$ be a nonempty subset of  a topological vector space $L$. We say that a set-valued function $F: X \to 2^L$ is a \say{Knaster-Kuratowski-Mazurkiewicz} {\bf $($KKM$)$} map if for each finite subset $\{x_1, x_2, \cdots, x_n\}$ of $X$ we have 
	\begin{equation}\label{KKM_inc} 
		{\rm co}(\{x_1, x_2, \cdots, x_n\}) \subset \bigcup\limits_{i = 1}^n F(x_i).
	\end{equation} 
\end{definition}

Notice that KKM maps are not necessarily lower semicontinuous even if $X$ is finite dimensional. Consider for example $F: [0,1] \to 2^{\mathbb{R}}$ given by $F(0) = [0,1]$ and $F(x) = [0,x]$ if $x \in (0,1].$ Then $F$ is closed and convex-valued, has a closed graph and is a KKM map. However, it fails to be lower semicontinuous. 

Moreover, there are KKM maps which are not graph-convex. For example $F:[0,1] \to 2^{\mathbb{R}}$ given by $F(x) = [0,x]$ if $x \in [0,1/2]$ and $F(x) = [x,1]$ if $x \in (1/2,1]$ is a KKM map which is not graph-convex.

Among the many properties of KKM maps (see for instance \cite{beer_topologies_1993} and references therein) is the fact that $x \in F(x)$ for all $x \in X$, i.e., $F$ is nonempty-valued and all vectors in $X$ are fixed points.

\begin{theorem}\label{kkm_theo}
	Assume that the function $\B{X}: X^\ad \to 2^{X^\ad}$ is a KKM map. Then the game generating such a map has an equilibrium.
\end{theorem}

\begin{proof}
	The proof is analogous to the proof of \emph{\Cref{no_lowersemicontinuity_theo_0}}. The only detail that changes is that before we used the graph-convexity of $\B{X}$ to show that convex combinations of fixed points were fixed points. In this case, since every point in $X^\ad$ is a fixed point, this step is not necessary.
\end{proof}

\section{Multiplier bias and uniqueness of an equilibrium for jointly convex games}\label{bias}

When each $\B{X}_i$ is given by \eqref{constraint_map} for some
$\C{C} \subset X$, the set of fixed points of $\B{X}$ coincides with $\C{C}$. This observation is important because every generalized Nash
 equilibrium must satisfy $\overline{x} \in \B{X}(\overline{x})$, and therefore
 $\overline{x} \in \C{C}$. In particular, if $\C{C}$ is closed and convex, the corresponding game is called \emph{jointly convex}. Jointly convex games allow for a notion of a subclass of Nash equilibria known as \emph{variational equilibria}. This notion is known to be more restrictive - in the sense that many more Nash equilibria may exist - but more amenable to numerical computations. Hence, it is very important for applied game theory \cite{facchinei_generalized_2007,hintermuller_generalized_2015,dvurechensky2024cournotnashmodelcoupledhydrogen}.

 \begin{definition}[\bf Variational equilibrium]
 	A strategy bundle $\overline{x} \in \C{C}$ is called a variational (or sometimes normalized) Nash equilibrium for a game $G = (\B{X}, \C{J})$ provided 
 	\begin{itemize}
 		\item[\bf (i)] it solves the minimization problem 
 		\begin{equation}\label{min_prob_eq_F} 
 			\min \Psi(\overline{x}, y) \ \ {\rm s.t. } \ \ y \in \C{C};
 		\end{equation} 
 	\end{itemize} or equivalently $($under $C^1$-regularity of the objectives$)$: 
 	\begin{itemize}
 		\item[\bf (ii)] the inequality 
 		\begin{equation}\label{qvi_game2n} 
 			\sum\limits_{i = 1}^N \langle \partial_i\C{J}_i(\overline{x}), y_i - \overline{x}_i \rangle_{X_i^*, X_i} \geqslant 0 
 		\end{equation} holds for all $y \in \C{C}.$
 	\end{itemize}
 \end{definition} 
\begin{remark}
 	For simplicity we have implicitly assumed that $\C{C} \subset X^\ad$, otherwise we need to replace $\C{C}$ with $X^\ad \cap \C{C}$ in the theorem above. 
\end{remark}
 Obviously variational equilibria are Nash equilibria. The converse is not true and there are multiple examples in the literature showing this fact, see for instance \cite{facchinei_generalized_2007,harker_generalized_1991}. 

In many applications, the set of equilibria of a GNEP is large and may contain continua of points \cite{facchinei_generalized_2007}. We want to end this section by extending Rosen’s diagonal strict convexity condition to our Banach space setting for jointly convex games. The main message is twofold:
(i) suitable structural conditions on a weighted pseudogradient single out a \emph{unique} variational equilibrium, and
(ii) changing the weights $r\in(\mathbb{R}_+^\ast)^N$ can be interpreted as introducing a \emph{multiplier bias} that selects one equilibrium among many. We keep the objective functionals fixed and bias only the way players react to marginal costs.

Let $G = (\B{X}, \C{J})$ be a jointly convex game in which the family of objective functionals satisfies \Cref{ass_objectives}. Consider the jointly convex game $G_r = (\B{X}, r \cdot \C{J})$ for some $r \in \mathbb{R}^N$. Here, '$\cdot$' is the componentwise multiplication of vectors. A well-known fact \cite{rosen_existence_1965} is that if $r \in (\mathbb{R}_+^*)^N$ (i.e., $r_i > 0$ for all $i \in I$), then  
\begin{equation}
	\B{E}(G_r) = \B{E}(G),
\end{equation}
where $\B{E}(\cdot)$ denotes the set of equilibria. This means that the set of Nash equilibria is unaffected by positive rescaling of the objective functions. An interesting aspect, however, is that except in the orthogonal case, i.e., when $\C{C}$ has the form $\C{C} = \C{C}_1 \times \cdots \times \C{C}_N$, we have $\B{VE}(G_r) \neq \B{VE}(G)$, where $\B{VE}(\cdot)$ is the set of variational equilibria. We show below, extending the results of \cite{rosen_existence_1965} to infinite dimensions, that an $r$-dependent structural condition on $G_r$ ensures that $\B{VE}(G_r)$ is a singleton. Hence, any numerical method designed to compute a variational equilibrium will converge to the same one. Since $\B{VE}(G_r) \subset \B{E}(G_r) = \B{E}(G)$ and $r \in (\mathbb{R}_+^*)^N$ must be chosen so that the structural condition holds, we are still computing \emph{a} Nash equilibrium of the original game.

Assume that the objective functionals $\{\C{J}_i\}_{i = 1}^N$ are continuously differentiable, with derivatives $\partial_i \C{J}_i(\cdot, x_\mt): X_i \to X_i^*$. Given a vector $r \in (\mathbb{R}_+^*)^N$, we define (formally) the \emph{pseudogradient} of $r \cdot \C{J}$ at $x \in X$ in the direction $h \in X$ as
\begin{equation}
	\label{sub_grad} d(x,r)h = \begin{bmatrix}
		r_1 \partial_1 \C{J}_1(x)h_1 \\ \vdots \\ r_N \partial_N\C{J}_N (x)h_N
	\end{bmatrix}.
\end{equation}
\begin{definition}
	We say that the game $G_r$ is diagonally strictly convex (with respect to the shared constraint $\C{C}$) for a given $r \in (\mathbb{R}_+^*)^N$ if for all $x, y \in \C{C}$,
	\begin{equation}
		d(x,r)(y-x) + d(y,r)(x-y) < 0.
	\end{equation}
\end{definition}

\begin{remark}
	In smoother problems, one can verify whether $G_r$ is diagonally strictly convex on $\C{C}$ by checking the positive definiteness of the operator $D(x,r) + D(x,r)^*$ for all $x \in \C{C}$, where $D(x,r)$ denotes the Jacobian of the map $x \mapsto d(x,r)$. Here, the superscript $^*$ indicates the adjoint operator.
\end{remark}

We begin with the non-generalized situation, i.e., we assume that $\B{X}(x) \equiv X^\ad$ even if $\B{X}_i$ has the form \eqref{constraint_map}. This corresponds to the case in which $\C{C}$ is a Cartesian product. The result below extends \cite[Theorem 2]{rosen_existence_1965} to the infinite-dimensional setting. In its proof we need the normal cone to $\C{C}$ at some $x\in\C{C}$ which is given by
$$
N_{\C{C}}(x)=\{\mu_x\in X^*:\mu_x(y-x)\leq 0\:\forall y\in\C{C}\}.
$$

\begin{prop}\label{propo_uniq}
	Assume that $\B{X}(x) \equiv X^\ad = \C{C}$ and that the family of objective functions satisfies \emph{\Cref{ass_objectives}}. If $G_r$ is diagonally strictly convex on $\C{C}$ for some $r \in \mathbb{R}_+^N$, then $\B{VE}(G)$ is a singleton.
\end{prop}

\begin{proof}
	From convex optimization along with the Nikaido-Isoda function formulation of a Nash equilibrium, it follows that if $\overline{x} \in \C{C}$ is a Nash equilibrium, then there exists $\mu_{\overline{x}} \in N_\C{C}(\overline{x})$ such that
	$$\nabla \Psi(\overline{x}, \overline{x}) + \mu_{\overline{x}} = 0\quad\text{in }X^*.$$
	Suppose that $\overline{x}, \overline{y} \in \C{C}$ are two distinct Nash equilibria. Testing the equation above by $r\cdot(\overline{y} - \overline{x})$ (in the $\overline{x}$ case) and by $r\cdot(\overline{x} - \overline{y})$ (in the $\overline{y}$ case), and adding the resulting equations, we obtain $\alpha + \beta = 0$, where
	$$\alpha = d(\overline{x}, r)(\overline{y} - \overline{x}) + d(\overline{y}, r)(\overline{x} - \overline{y}),$$ $$
	\beta = \mu_{\overline{x}}(r\cdot(\overline{y} - \overline{x})) + \mu_{\overline{y}}(r\cdot(\overline{x} - \overline{y})).$$
	Since $r \in \mathbb{R}_+^N$, it follows from the definition of the normal cone and from the fact that $X^\ad = \C{C}$ that $\beta \leqslant  0$, and hence $\alpha \geqslant  0$, which contradicts diagonal strict convexity. Therefore, $\overline{x} = \overline{y}$.
\end{proof}

We now return to the case in which $\B{X}$ is non-constant but still given by \eqref{constraint_map}.  
\begin{prop}\label{propo_uniq2}
	Assume that $\B{X}$ has the structure \eqref{constraint_map} and that the family of objective functions satisfies \emph{\Cref{ass_objectives}}. If $G_r$ is diagonally strictly convex on $\C{C}$ for all $r \in \mathcal{R}$, where $\mathcal{R} \subset (\mathbb{R}_+^*)^N$ is nonempty and convex, then for each $r \in \mathcal{R}$ the set $\B{VE}(G_r) \subset \B{E}(G)$ is a singleton.
\end{prop}

It follows from \Cref{propo_uniq} and \Cref{propo_uniq2} that when the constraint map is constant, the variational equilibrium is unique and independent of $r$. However, in the non-constant case, although the scaling vector $r$ cannot alter the set of Nash equilibria, it may affect the set of variational equilibria – a feature of practical interest given their computational advantages over the full equilibrium set. The next proposition quantifies this effect.

\begin{prop}\label{propo_dir}
	Suppose that $G_r$ is diagonally strictly convex for all $r \in \mathcal{R}.$ Let $r, s \in \mathcal{R}\subset (\mathbb{R}_+^*)^N$ satisfy $r_j > s_j$ for some $j \in I$ and $r_i = s_i$ for all $i \neq j$, and let $x^r, x^s$ denote the corresponding variational equilibria. Then the directional derivative of $\C{J}_{j}$ at $x^r$ in the direction $x_{j}^r - x_{j}^s$ is negative.
\end{prop}

We do not include the proofs of \Cref{propo_uniq2} and \Cref{propo_dir} because the arguments follow the same structure as \Cref{propo_uniq}. The corresponding finite-dimensional proofs are found in \cite{rosen_existence_1965}.

Diagonal strict convexity also implies a strong stability property, naturally formulated in terms of an associated ordinary differential equation type dynamical system that generates a trajectory of strategies. This dynamical perspective has served as a foundation for algorithms computing variational equilibria and offers a complementary viewpoint to modern numerical methods for variational inequalities; see, for example, \cite{facchinei_generalized_2007,facchinei_decomposition_2011} and references therein. Beyond its algorithmic relevance, this interpretation provides a useful conceptual and historical connection between equilibrium theory and dynamical systems.

\section{Conclusions and Outlook}

This work revisited the existence theory for convex generalized Nash equilibrium problems in Banach spaces. We showed that equilibrium existence does not require lower semicontinuity of the constraint maps, a classical but restrictive assumption, and that graph-convexity or the KKM property are sufficient. These conditions are geometric in nature and easier to verify in PDE-constrained and infinite-dimensional games.

Future research may address quantitative stability of equilibria under perturbations of the constraint maps, algorithmic schemes exploiting graph-convex or KKM structures, and extensions to nonconvex or stochastic settings. An additional direction suggested by the geometric formulation is to investigate mixed analytic–geometric models in which some players are described by objective functionals and others by preference correspondences, a situation that naturally arises in multi-agent systems with heterogeneous information or incomplete preference specification.

\bibliographystyle{abbrv} 
\bibliography{references}

@book {MR895589,
    AUTHOR = {Wloka, J.},
     TITLE = {Partial differential equations},
      NOTE = {Translated from the German by C. B. Thomas and M. J. Thomas},
 PUBLISHER = {Cambridge University Press, Cambridge},
      YEAR = {1987},
     PAGES = {xii+518},
      ISBN = {0-521-25914-2; 0-521-27759-0},
   MRCLASS = {35-01 (46-01 65-01)},
  MRNUMBER = {895589},
       DOI = {10.1017/CBO9781139171755},
       URL = {https://doi.org/10.1017/CBO9781139171755},
}

@book {MR2424078,
    AUTHOR = {Adams, Robert A. and Fournier, John J. F.},
     TITLE = {Sobolev spaces},
    SERIES = {Pure and Applied Mathematics (Amsterdam)},
    VOLUME = {140},
   EDITION = {Second},
 PUBLISHER = {Elsevier/Academic Press, Amsterdam},
      YEAR = {2003},
     PAGES = {xiv+305},
      ISBN = {0-12-044143-8},
   MRCLASS = {46E35 (46-01 46-02 46B70 46Exx)},
  MRNUMBER = {2424078},
}

@incollection{hintermuller_recent_2019,
	address = {Cham},
	title = {Recent {Trends} and {Views} on {Elliptic} {Quasi}-{Variational} {Inequalities}},
	isbn = {978-3-030-33115-3 978-3-030-33116-0},
	url = {http://link.springer.com/10.1007/978-3-030-33116-0_1},
	language = {en},
	urldate = {2025-07-03},
	booktitle = {Topics in {Applied} {Analysis} and {Optimisation}},
	publisher = {Springer International Publishing},
	author = {Alphonse, Amal and Hintermüller, Michael and Rautenberg, Carlos N.},
	editor = {Hintermüller, Michael and Rodrigues, José Francisco},
	year = {2019},
	doi = {10.1007/978-3-030-33116-0_1},
	note = {Series Title: CIM Series in Mathematical Sciences},
	pages = {1--31},
}

@article{arrow_existence_1954,
	title = {Existence of an {Equilibrium} for a {Competitive} {Economy}},
	volume = {22},
	issn = {0012-9682},
	url = {https://www.jstor.org/stable/1907353},
	doi = {10.2307/1907353},
	number = {3},
	urldate = {2024-05-21},
	journal = {Econometrica},
	author = {Arrow, Kenneth J. and Debreu, Gerard},
	year = {1954},
	note = {Publisher: [Wiley, Econometric Society]},
	pages = {265--290},
}

@book{beer_topologies_1993,
	address = {Dordrecht},
	edition = {1st ed. 1993},
	series = {Mathematics and {Its} {Applications}},
	title = {Topologies on {Closed} and {Closed} {Convex} {Sets}},
	isbn = {978-90-481-4333-7 978-94-015-8149-3},
	publisher = {Springer Netherlands},
	author = {Beer, Gerald},
	year = {1993},
	doi = {10.1007/978-94-015-8149-3},
}

@book{bensoussan_impulse_1987,
	series = {Modern {Applied} {Mathematics} {Series}},
	title = {Impulse {Control} and {Quasi} {Variational} {Inequalities}},
	isbn = {978-0-471-82973-7},
	publisher = {John Wiley \& Sons Canada, Limited},
	author = {Bensoussan, Alain and Lions, Jacques-Louis},
	year = {1987},
}

@article{MR3032892,
    AUTHOR = {Borz\`i, Alfio and Kanzow, Christian},
     TITLE = {Formulation and numerical solution of {N}ash equilibrium
              multiobjective elliptic control problems},
   JOURNAL = {SIAM J. Control Optim.},
  FJOURNAL = {SIAM Journal on Control and Optimization},
    VOLUME = {51},
      YEAR = {2013},
    NUMBER = {1},
     PAGES = {718--744},
      ISSN = {0363-0129,1095-7138},
   MRCLASS = {49N70 (49K20 49M25 65K10 91A23)},
  MRNUMBER = {3032892},
MRREVIEWER = {Philip\ D.\ Loewen},
       DOI = {10.1137/120864921},
       URL = {https://doi.org/10.1137/120864921},
}

@article{breton_game-theoretic_2006,
	title = {A game-theoretic formulation of joint implementation of environmental projects},
	volume = {168},
	issn = {0377-2217},
	url = {https://www.sciencedirect.com/science/article/pii/S0377221704007248},
	doi = {10.1016/j.ejor.2004.06.045},
	number = {1},
	urldate = {2024-07-19},
	journal = {European Journal of Operational Research},
	author = {Breton, Michèle and Zaccour, Georges and Zahaf, Mourad},
	month = jan,
	year = {2006},
	pages = {221--239},
}

@article{chan_generalized_1982,
	title = {The {Generalized} {Quasi}-{Variational} {Inequality} {Problem}},
	volume = {7},
	issn = {0364-765X, 1526-5471},
	url = {https://pubsonline.informs.org/doi/10.1287/moor.7.2.211},
	doi = {10.1287/moor.7.2.211},
	number = {2},
	urldate = {2024-07-19},
	journal = {Mathematics of Operations Research},
	author = {Chan, Der-San and Pang, Jong-Shi},
	month = may,
	year = {1982},
	pages = {211--222},
}

@misc{dvurechensky2024cournotnashmodelcoupledhydrogen,
  title        = {A Cournot-Nash Model for a Coupled Hydrogen and Electricity Market},
  author       = {Dvurechensky, Pavel and Geiersbach, Caroline and Hintermüller, Michael and Kannan, Aswin and Kater, Stefan and Zöttl, Gregor},
  year         = {2024},
  eprint       = {2410.20534},
  archivePrefix= {arXiv},
  primaryClass = {math.OC},
  url          = {https://arxiv.org/abs/2410.20534},
  note         = {arXiv:2410.20534 [math.OC]},
}

@article{facchinei_generalized_2007,
	title = {Generalized {Nash} equilibrium problems},
	volume = {5},
	issn = {1619-4500},
	url = {https://doi.org/10.1007/s10288-006-0016-9},
	doi = {10.1007/s10288-006-0016-9},
	number = {3},
	urldate = {2024-05-21},
	journal = {4OR},
	author = {Facchinei, Francisco and Kanzow, Christian},
	month = sep,
	year = {2007},
	pages = {173--210},
}

@article{facchinei_decomposition_2011,
	title = {Decomposition algorithms for generalized potential games},
	volume = {50},
	issn = {0926-6003, 1573-2894},
	url = {https://doi.org/10.1007/s10589-010-9357-7},
	doi = {10.1007/s10589-010-9357-7},
	number = {2},
	urldate = {2024-05-21},
	journal = {Computational Optimization and Applications},
	author = {Facchinei, Francisco and Piccialli, Veronica and Sciandrone, Marco},
	month = oct,
	year = {2011},
	pages = {237--262},
}

@article{fan_minimax_1953,
	title = {Minimax {Theorems}},
	volume = {39},
	issn = {0027-8424, 1091-6490},
	url = {https://www.pnas.org/doi/10.1073/pnas.39.1.42},
	doi = {10.1073/pnas.39.1.42},
	number = {1},
	urldate = {2024-07-19},
	journal = {Proceedings of the National Academy of Sciences},
	author = {Fan, Ky},
	month = jan,
	year = {1953},
	pages = {42--47},
}

@incollection{gahururu_generalized_2022,
	address = {Cham},
	title = {Generalized {Nash} {Equilibrium} {Problems} with {Partial} {Differential} {Operators}: {Theory}, {Algorithms}, and {Risk} {Aversion}},
	isbn = {978-3-030-82184-5 978-3-030-82185-2},
	url = {http://link.springer.com/10.1007/978-3-030-82185-2_6},
	language = {en},
	urldate = {2025-07-03},
	booktitle = {Non-{Smooth} and {Complementarity}-{Based} {Distributed} {Parameter} {Systems}: {Simulation} and {Hierarchical} {Optimization}},
	publisher = {Springer International Publishing},
	author = {Gahururu, Didier and Hintermüller, Michael and Stengl, Steven-Marian and Surowiec, Thomas M.},
	editor = {Hintermüller, Michael and Herzog, Ren{\'e} and Kanzow, Christian and Ulbrich, Michael and Ulbrich, Stefan},
	year = {2022},
	doi = {10.1007/978-3-030-82185-2_6},
	pages = {145--181},
}

@article{gahururu_risk-neutral_2023,
	title = {Risk-neutral {PDE}-constrained generalized {Nash} equilibrium problems},
	volume = {198},
	issn = {0025-5610, 1436-4646},
	url = {https://link.springer.com/10.1007/s10107-022-01916-2},
	doi = {10.1007/s10107-022-01916-2},
	number = {2},
	urldate = {2025-07-03},
	journal = {Mathematical Programming},
	author = {Gahururu, Didier B. and Hintermüller, Michael and Surowiec, Thomas M.},
	month = apr,
	year = {2023},
	pages = {1287--1337},
}

@article{harker_generalized_1991,
	title = {Generalized {Nash} games and quasi-variational inequalities},
	volume = {54},
	issn = {0377-2217},
	url = {https://www.sciencedirect.com/science/article/pii/037722179190205E},
	doi = {10.1016/0377-2217(91)90205-E},
	number = {1},
	urldate = {2024-07-19},
	journal = {European Journal of Operational Research},
	author = {Harker, Patrick T.},
	month = sep,
	year = {1991},
	pages = {81--94},
}

@book{harsanyi_general_1988,
	series = {MIT {Press} {Books}},
	title = {A {General} {Theory} of {Equilibrium} {Selection} in {Games}},
	number = {1},
	publisher = {The MIT Press},
	author = {Harsanyi, John C. and Selten, Reinhard},
	year = {1988},
}

@article{hintermueller_pde-constrained_2013,
	title = {A {PDE}-{Constrained} {Generalized} {Nash} {Equilibrium} {Problem} with {Pointwise} {Control} and {State} {Constraints}},
	volume = {9},
	issn = {1348-9151},
	number = {2},
	journal = {Pacific Journal of Optimization},
	author = {Hinterm\"uller, Michael and Surowiec, Thomas},
	month = apr,
	year = {2013},
	pages = {251--273},
}

@article{hintermuller_generalized_2024,
  title = {A Generalized {$\Gamma$}-Convergence Concept for a Class of Equilibrium Problems},
  volume = {34},
  issn = {1432-1467},
  url = {https://doi.org/10.1007/s00332-024-10059-x},
  doi = {10.1007/s00332-024-10059-x},
  language = {en},
  number = {5},
  urldate = {2025-08-05},
  journal = {Journal of Nonlinear Science},
  author = {Hintermüller, Michael and Stengl, Steven-Marian},
  month = jul,
  year = {2024},
  pages = {83},
}

@article{hintermuller_generalized_2015,
	title = {Generalized {Nash} {Equilibrium} {Problems} in {Banach} {Spaces}: {Theory}, {Nikaido}--{Isoda}-{Based} {Path}-{Following} {Methods}, and {Applications}},
	volume = {25},
	issn = {1052-6234, 1095-7189},
	url = {http://epubs.siam.org/doi/10.1137/130929203},
	doi = {10.1137/130929203},
	number = {3},
	urldate = {2025-07-03},
	journal = {SIAM Journal on Optimization},
	author = {Hintermüller, Michael and Surowiec, Thomas and Kämmler, Alexander},
	month = jan,
	year = {2015},
	pages = {1826--1856},
}

@article{hobbs_nash-cournot_2007,
	title = {Nash-{Cournot} {Equilibria} in {Electric} {Power} {Markets} with {Piecewise} {Linear} {Demand} {Functions} and {Joint} {Constraints}},
	volume = {55},
	issn = {0030-364X, 1526-5463},
	url = {https://pubsonline.informs.org/doi/10.1287/opre.1060.0315},
	doi = {10.1287/opre.1060.0315},
	number = {1},
	urldate = {2024-07-19},
	journal = {Operations Research},
	author = {Hobbs, Benjamin F. and Pang, Jong-Shi},
	month = feb,
	year = {2007},
	pages = {113--127},
}

@article{jing-yuan_spatial_1999,
	title = {Spatial {Oligopolistic} {Electricity} {Models} with {Cournot} {Generators} and {Regulated} {Transmission} {Prices}},
	volume = {47},
	issn = {0030-364X, 1526-5463},
	url = {https://pubsonline.informs.org/doi/10.1287/opre.47.1.102},
	doi = {10.1287/opre.47.1.102},
	number = {1},
	urldate = {2024-07-19},
	journal = {Operations Research},
	author = {Jing-Yuan, Wang and Smeers, Yves},
	month = feb,
	year = {1999},
	pages = {102--112},
}

@article{kanzow_multiplier-penalty_2019,
	title = {The {Multiplier}-{Penalty} {Method} for {Generalized} {Nash} {Equilibrium} {Problems} in {Banach} {Spaces}},
	volume = {29},
	issn = {1052-6234, 1095-7189},
	url = {http://epubs.siam.org/doi/10.1137/17M1124735},
	doi = {10.1137/17M1124735},
	number = {1},
	urldate = {2024-05-21},
	journal = {SIAM Journal on Optimization},
	author = {Kanzow, Christian and Karl, Victor and Steck, Daniel and Wachsmuth, Daniel},
	month = jan,
	year = {2019},
	pages = {767--793},
}

@incollection{mas-colell_equilibrium_1991,
	series = {Handbooks in {Econom}.},
	title = {Equilibrium theory in infinite dimensional spaces},
	volume = {4},
	booktitle = {Handbook of {Mathematical} {Economics}},
	publisher = {North-Holland},
	author = {Mas-Colell, Andreu and Zame, William R.},
	year = {1991},
}

@article{nash_non-cooperative_1951,
	title = {Non-{Cooperative} {Games}},
	volume = {54},
	issn = {0003-486X, 1939-8980},
	url = {https://www.jstor.org/stable/1969529},
	doi = {10.2307/1969529},
	number = {2},
	urldate = {2024-07-19},
	journal = {Annals of Mathematics},
	author = {Nash, John},
	month = sep,
	year = {1951},
	pages = {286},
}

@article{nash_equilibrium_1950,
	title = {Equilibrium points in n-person games},
	volume = {36},
	issn = {0027-8424, 1091-6490},
	url = {https://www.pnas.org/doi/10.1073/pnas.36.1.48},
	doi = {10.1073/pnas.36.1.48},
	number = {1},
	urldate = {2024-07-19},
	journal = {Proceedings of the National Academy of Sciences},
	author = {Nash, John F.},
	month = jan,
	year = {1950},
	pages = {48--49},
}

@article{nikaido_note_1955,
	title = {Note on non-cooperative convex game},
	volume = {5},
	issn = {0030-8730},
	url = {https://projecteuclid.org/journals/pacific-journal-of-mathematics/volume-5/issue-5/Note-on-non-cooperative-convex-game/pjm/1103045138.full},
	doi = {10.2140/pjm.1955.5.807},
	number = {5},
	urldate = {2024-07-19},
	journal = {Pacific Journal of Mathematics},
	author = {Nikaidô, Hukukane and Isoda, Kazuo},
	month = dec,
	year = {1955},
	pages = {807--815},
}

@article{pang_quasi-variational_2005,
	title = {Quasi-variational inequalities, generalized {Nash} equilibria, and multi-leader-follower games},
	volume = {2},
	issn = {1617-5891, 1617-5905},
	url = {https://doi.org/10.1007/s10287-004-0020-y},
	doi = {10.1007/s10287-004-0020-y},
	number = {1},
	urldate = {2024-05-21},
	journal = {Computational Management Science},
	author = {Pang, Jong-Shi and Fukushima, Masao},
	month = jan,
	year = {2005},
	pages = {21--56},
}

@article{ramos_nash_2002,
	title = {Nash {Equilibria} for the {Multiobjective} {Control} of {Linear} {Partial} {Differential} {Equations}},
	volume = {112},
	issn = {0022-3239, 1573-2878},
	url = {https://doi.org/10.1023/A:1015654401327},
	doi = {10.1023/A:1015654401327},
	number = {3},
	urldate = {2024-07-19},
	journal = {Journal of Optimization Theory and Applications},
	author = {Ramos, A. and Glowinski, Roland and Periaux, Jacques},
	month = mar,
	year = {2002},
	pages = {457--498},
}

@article{ramos_pointwise_2002,
	title = {Pointwise {Control} of the {Burgers} {Equation} and {Related} {Nash} {Equilibrium} {Problems}: {Computational} {Approach}},
	volume = {112},
	issn = {0022-3239, 1573-2878},
	url = {https://doi.org/10.1023/A:1015667120037},
	doi = {10.1023/A:1015667120037},
	number = {3},
	urldate = {2024-07-19},
	journal = {Journal of Optimization Theory and Applications},
	author = {Ramos, A. and Glowinski, Roland and Periaux, Jacques},
	month = mar,
	year = {2002},
	pages = {499--516},
}

@book{rockafellar_variational_1998,
	series = {Grundlehren der mathematischen {Wissenschaften}},
	title = {Variational {Analysis}},
	isbn = {978-3-540-62772-2 978-3-642-02431-3},
	shorttitle = {Variational {Analysis}},
	publisher = {Springer Berlin Heidelberg},
	author = {Rockafellar, R. Tyrrell and Wets, Roger J. B.},
	year = {1998},
	doi = {10.1007/978-3-642-02431-3},
	note = {ISSN: 0072-7830},
}

@article{rosen_existence_1965,
	title = {Existence and {Uniqueness} of {Equilibrium} {Points} for {Concave} {N}-{Person} {Games}},
	volume = {33},
	issn = {0012-9682},
	url = {https://www.jstor.org/stable/1911749},
	doi = {10.2307/1911749},
	number = {3},
	urldate = {2024-07-19},
	journal = {Econometrica},
	author = {Rosen, J. B.},
	year = {1965},
	pages = {520--534},
}

@incollection{roubicek_noncooperative_1999,
	address = {Basel},
	series = {International {Series} of {Numerical} {Mathematics}},
	title = {Noncooperative {Games} with {Elliptic} {Systems}},
	isbn = {978-3-0348-8683-3 978-3-0348-8690-1},
	url = {http://link.springer.com/10.1007/978-3-0348-8690-1_14},
	language = {en},
	urldate = {2024-05-21},
	booktitle = {Optimal {Control} of {Partial} {Differential} {Equations}},
	publisher = {Birkhäuser},
	author = {Roubíček, Tomáš},
	editor = {Hoffmann, Karl-Heinz and Leugering, Günter and Tröltzsch, Fredi and Caesar, Siegfried},
	year = {1999},
	doi = {10.1007/978-3-0348-8690-1_14},
	pages = {245--255},
}

@article{roubicek_nash_2007,
	title = {On {Nash} {Equilibria} for {Noncooperative} {Games} {Governed} by the {Burgers} {Equation}},
	volume = {132},
	issn = {0022-3239, 1573-2878},
	url = {https://doi.org/10.1007/s10957-007-9238-y},
	doi = {10.1007/s10957-007-9238-y},
	number = {1},
	urldate = {2024-07-19},
	journal = {Journal of Optimization Theory and Applications},
	author = {Roubíček, Tomáš},
	month = mar,
	year = {2007},
	pages = {41--50},
}

@article{roughgarden_algorithmic_2010,
	title = {Algorithmic game theory},
	volume = {53},
	issn = {0001-0782, 1557-7317},
	url = {https://dl.acm.org/doi/10.1145/1785414.1785444},
	doi = {10.1145/1785414.1785444},
	number = {7},
	urldate = {2024-07-19},
	journal = {Communications of the ACM},
	author = {Roughgarden, Tim},
	month = jul,
	year = {2010},
	pages = {78--86},
}

@book{von_neumann_theory_2007,
	series = {A {Princeton} classic edition},
	title = {Theory of games and economic behavior},
	isbn = {978-0-691-13002-4},
	publisher = {Princeton University Press},
	address = {Princeton, NJ},
	edition = {60. anniversary ed., 4. print., and 1. paperb. print},
	author = {Von Neumann, John and Morgenstern, Oskar},
	year = {2007},
}

@article{bongarti_optimal_2024,
    title = {Optimal {Boundary} {Control} of the {Isothermal} {Semilinear} {Euler} {Equation} for {Gas} {Dynamics} on a {Network}},
    volume = {89},
    copyright = {All rights reserved},
    issn = {0095-4616, 1432-0606},
    url = {https://link.springer.com/10.1007/s00245-023-10088-0},
    doi = {10.1007/s00245-023-10088-0},
    language = {en},
    number = {2},
    urldate = {2024-05-22},
    journal = {Applied Mathematics \& Optimization},
    author = {Bongarti, Marcelo and Hintermüller, Michael},
    month = apr,
    year = {2024},
    pages = {36},
}
\end{document}